\documentclass[12pt]{amsart}
\usepackage{mathptmx}

\let\oldlabel=\label
\def\prellabel{\def\label##1{\oldlabel{##1}\ifmmode\else\ifinner\else
          \marginpar{{\footnotesize\ \\ \tt
                     ##1}}\fi\fi}}
\def\discuss#1{\relax}

\def\CHANGE{\discuss{} }
\def\discuss#1{\marginpar{\raggedright\tiny #1}}

\def\NN{{\NZQ N }}

\def\ZZ{{\NZQ Z}}
\def\RR{{\NZQ R}}

\def\opn#1#2{\def#1{\operatorname{#2}}} % to make operators
\opn\chara{char}
\opn\gr{gr} \opn\rank{rank}
\let\union=\cup
\let\sect=\cap

\let\tensor=\otimes
\let\iso=\cong

\let\Dirsum=\bigoplus

\let\:=\colon
\newtheorem{theorem}{Theorem}[section]
\newtheorem{lemma}[theorem]{Lemma}
\newtheorem{corollary}[theorem]{Corollary}
\newtheorem{proposition}[theorem]{Proposition}
\theoremstyle{definition}
\newtheorem{remark}[theorem]{Remark}

\newtheorem{definition}[theorem]{Definition}
\let\epsilon=\varepsilon
\let\phi=\varphi
\opn\ini{in} \opn\KRS{KRS}

\def\pp{{\mathfrak p}}
\def\qq{{\mathfrak q}}
\opn\krs{krs}
\opn\diag{diag}
\opn\DD{{\mathcal D}}
\opn\SS{{\mathcal S}}
\opn\MM{{\mathcal M}}
\opn\GL{GL}

\let\Bbb=\mathbb

\def\RR{{\Bbb R}}

\def\ZZ{{\Bbb Z}}
\def\NN{{\Bbb N}}
\def\ini{\operatorname{in}}

\opn\height{height}
\opn\length{length}
\opn\cl{cl}
\opn\Cl{Cl}
\opn\Grass{Grass}
\opn\sr{sr}
\opn\Hom{Hom}
\def\sep{\,|\,}
\opn\Ker{Ker}
\opn\Im{Im}

\opn\supp{supp}
\def\ol#1{\overline{#1}}
\opn\Spec{Spec}
\def\OO{{\mathcal O}}

\opn\SL{SL}
\opn\Sing{Sing}

\opn\QF{QF}

\def\LL{{\mathcal L}}

\def\FF{{\mathcal F}}
\def\BB{{\mathcal B}}
\def\dd{{\mathbf d}}

\def\bl{\bigl[}
\def\br{\bigr]}
\def\bs{\ \big|\ }

\def\hak{\mathbin{\unitlength=1mm\thicklines\,\begin{picture}%
     (2,2)(0,0)\put(0,0.3){\line(0,1){2.1}}\put(0,0.3){\line(1,0){1.8}}\end{picture}\,}}

\begin{document}

\title{The variety of exterior powers of linear maps}
\author{Winfried Bruns and Aldo Conca}
\address{Universit\"at Osnabr\"uck, FB
Mathematik/Informatik, 49069 Osna\-br\"uck, Germany}
\email{wbruns@uos.de}
\address{Dipartimento di Matematica, Universit\'a di Genova,
Via Dodeca\-neso 35, 16146 Genova, Italy}
\email{conca@dima.unige.it}

\begin{abstract}
Let $V$ and $W$ be vector spaces of dimension $m$ and $n$ resp. We
investigate the Zariski closure $X_t$ of the image $Y_t$ of the map
$\Hom_K(V,W)\to\Hom_K(\bigwedge^t V,\bigwedge^t W)$, $\phi\mapsto
\bigwedge^t \phi$. In the case $t=\min(m,n)$, $Y_t=X_t$ is the cone
over a Grassmannian, but for $1<t<\min(m,n)$ one has $X_t\neq Y_t$.
We analyze the $G=\GL(V)\times\GL(W)$-orbits in $X_t$ via the
$G$-stable prime ideals in $\OO(X_t)$. It turns out that they are
classified by two numerical invariants, one of which is the rank and
the other a related invariant that we call small rank. Surprisingly,
the orbits in $X_t\setminus Y_t$ arise from the images $Y_u$ for
$u<t$ and simple algebraic operations. In the last section we
determine the singular locus of $X_t$. Apart from well-understood
exceptional cases, it is formed by the elements of rank $\le 1$ in
$Y_t$.
\end{abstract}

\maketitle

\section{Introduction}

Let $K$ be a field, and $V,W$ vector spaces over $K$ of finite
dimensions $m$ and $n$. In this paper we want to study the exterior
power map
$$
\Lambda_t: \Hom_K(V,W) \to \Hom_K\biggl(\bigwedge^t V, \bigwedge^t
W\biggr),\qquad \Lambda_t(\alpha)=\bigwedge^t \alpha.
$$
We want to analyze the Zariski closure $X_t(V,W)$ of the image
$Y_t(V,W)$ of $\Lambda_t$. In the following we will always assume
that $m\le n$. This does not restrict the generality since
$\Lambda_t$ commutes with dualization.

There are three special cases that make it impossible to formulate
all results in a uniform way. In the first two of them,
$X_t(V,W)=\Hom_K(\bigwedge^t V, \bigwedge^t W)$ so that $X_t(V,W)$
is an affine space: (i) in the trivial case $t=1$, and (ii) in the
case $t=m-1=n-1$; in fact, one has $\dim X_t(V,W)=mn$, unless
$t=m>1$ (see \cite[(10.16)(b)]{BV} or Proposition \ref{image1}).
Therefore, if $t=m-1=n-1$, then
$$
\dim X_t(V,W)=m^2=\binom mt\binom mt=\dim \Hom_K\biggl(\bigwedge^t
V, \bigwedge^t W\biggr).
$$
The third case is highly nontrivial, but classical: (iii) if $t=m$
(and $m\le n$), then $Y_t(V,W)$ is the cone over the Grassmannian
$\Grass(t,W)$, and, in particular, it is Zariski closed (for example,
see \cite{BV}). In all cases different from (i) and (iii),
$Y_t(V,W)$ is strictly contained in $X_t(V,W)$, as we will see.

For a compact formulation of our main result let us identify
$\Hom_K(V,W)$ with $V^*\tensor W$ and $\Hom_K\bigl(\bigwedge^t V,
\bigwedge^t W\bigr)$ with $\bigwedge^t V^*\tensor \bigwedge^t W$.
Moreover we consider all $t$ simultaneously by taking the unions
$$
Y(V,W)=\bigcup_{t\ge0} Y_t(V,W)\quad\text{and}\quad
X(V,W)=\bigcup_{t\ge0} X_t(V,W)
$$
in the algebra $\bigwedge V^*\tensor \bigwedge W$.

\begin{theorem}\label{main1}
Let $K$ be an algebraically closed field of characteristic $0$. Then
$X(V,W)$ is the closure of $Y(V,W)$ under the operation of
$V^*\times W$ by multiplication on $\bigwedge V^*\tensor \bigwedge
W$.
\end{theorem}

Clearly, in the algebra $\bigwedge V^*\tensor \bigwedge W$ (or its
subalgebra $\Dirsum_t \bigwedge^t V^*\tensor \bigwedge^t W$) the map
$\Lambda_t$ is just the $t$th power map, but usually we prefer the
viewpoint of linear maps and their exterior powers.

The group $G=\GL(V)\times\GL(W)$ acts naturally on $Y_t(V,W)$ and
$X_t(V,W)$. The proof of Theorem \ref{main1} is based on an analysis
of the orbit structure of $X_t(V,W)$. It turns out that these orbits
are classified by two numerical invariants. One of them is the
ordinary rank of an element $x\in \Hom_K\bigl(\bigwedge^t V,
\bigwedge^t W\bigr)$. The other one is a new invariant that we will
call \emph{small rank} $\sr(x)$. \CHANGE As we will see, $\sr(x)$ indicates
from which $Y_u(V,W)$ the element $x$  ``originates"  in the
sense of Theorem \ref{main1}: for $x\in X_t(V,W)$ with $\rank x>1$ a
representation
$$
x=\bigl((\alpha_1\wedge\dots\wedge\alpha_{t-s})\tensor
(y_1\wedge\dots\wedge y_{t-s})\bigr)\cdot x'
$$
with $x'\in Y_s(V,W)$, $\alpha_1,\dots,\alpha_{t-s}\in V^*$ and
$y_1,\dots,y_{t-s}\in W$ is possible if and only if $s=\sr(x)-1$,
and if $\rank x=1$, then $\sr(x)=1$ as well, and one can choose
$S=0$.

The determination of the orbits is based on the construction of
normal forms for suitable pairs of rank and small rank on one side
(Section \ref{numinvs}), and an analysis of the coordinate ring
$A_t(V,W)$ of $X_t(V,W)$ on the other (Section \ref{Atstruct}).
After the choice of bases in $V$ and $W$, one can identify the
coordinate ring of $\Hom_K(V,W)$ with the polynomial ring $K[X]$ in
the entries of an $m\times n$ matrix $X$ of indeterminates over $K$
and $A_t(V,W)$ with the subalgebra $A_t=A_t(m,n)$ of $K[X]$
generated by the $t$-minors of $X$. Using the decomposition of $A_t$
into irreducible $G$-submodules, we will determine the $G$-stable
prime ideals in $A_t$ (Theorem \ref{all-primes}). At this point, the
hypothesis of characteristic $0$ is used in a crucial way.

We will then analyze the fibers of $\Lambda_t$ (Section
\ref{Lambda-t}). The normal forms make it easy to understand the
effect of the (iterated) multiplication by elements of $V^*\times
W$. It is then not difficult to determine the structure of the
orbits (Section \ref{orbits}) and to prove Theorem \ref{main1}. It
is an important point that $X_{t'}(V',W')$ is a retract of
$X_t(V,W)$ if $t-t'\ge \dim V-\dim V',\dim W-\dim W'$.

In the last part of the paper (Section \ref{Sing}) we determine the
singular locus of $X_t$. Apart from the special cases (i), (ii) and
(iii) described above, the singular locus of $X_t$ is formed by the
elements of rank $\le 1$. The singular locus is always contained in
$Y_t$.

Throughout the paper, $K$ is an algebraically closed field of
characteristic $0$. However, among the basic arguments, only the
determination of the $G$-stable prime ideals in $A_t(m,n)$ depends
on it. We believe that all the results can be extended without
changes to the case $\chara K>\min(t,m-t,n-t)$, which we call
\emph{non-exceptional}. The extension to arbitrary characteristic
may require some changes.

\medskip

\noindent\emph{Conventions.}\enspace For a more compact notation we
set $L=L(V,W)=\Hom_K(V,W)$ and
$$
\LL_t=\LL_t(V,W)=\Hom_K\biggl(\bigwedge^t V, \bigwedge^t W\biggr).
$$
It is clear that the dimensions $m$ and $n$ of $V$ and $W$ define
all our data up to isomorphism, and therefore we will often replace
$V$ and $W$ by them. For example we will write $X_t(m,n)$ for
$X_t(V,W)$, and often $X_t$ and $Y_t$ will denote our objects
unambiguously.

Similarly we will write $A_t(m,n)$ or $A_t$ for the coordinate ring
of $X_t$. The minors generating it are homogeneous elements of
degree $t$. We can therefore normalize degrees in $A_t$, dividing
the degree in $K[X]$ by $t$.

The group $G=\GL(V)\times\GL(W)$ acts naturally on $L$ via
$$
\beta\mapsto \gamma\circ\beta\circ\alpha^{-1},\qquad
(\alpha,\gamma)\in G.
$$
By the functoriality of the $t$-th exterior power it acts likewise
on $\LL$, and the map $\Lambda$ is $G$-equivariant. Consequently $G$
acts on $Y=Y_t(V,W)$ and on $X=X_t(V,W)$.

It will sometimes be useful to allow that $t=0$. By definition,
$\bigwedge^0 \beta$ is the identity on $K=\bigwedge^0 V=\bigwedge^0
W$ for all $\beta\in L$.

With a basis $e_1,\dots,e_m$ of $V$ we associate the basis of
$\bigwedge ^t V$ that consists of the products
$e_{i_1}\wedge\dots\wedge e_{i_t}$ with $i_1 < \dots < i_t$. If
necessary we order these basis elements lexicographically by their
indices. The same convention applies to $W$.

If an element $\beta\in L$ is represented by a matrix $B$ with
respect to given bases of $V$ and $W$, then, with the choice of
bases of the exterior powers just specified, $\bigwedge^t \beta$ is
given by the matrix $\bigwedge^t B$ whose entries are the $t$-minors
$$
[i_1\dots i_t\sep j_1\dots j_t]_B
$$
of $B$. This notation will be used for $t$-minors of matrices in
general, and the index $B$ may be omitted if no confusion arises.

\section{Numerical invariants and normal forms}\label{numinvs}
It is easy to see that $Y_t(m,n)$ consists of exactly $m-t+2$
orbits. In fact, $L$ consists of $m+1$ orbits characterized by the
ranks of the elements in $L$. Of these orbits, $t$ are mapped to
$0\in\LL$, namely those corresponding to the ranks $0,\dots,t-1$,
and the images of the remaining ones stay disjoint in $\LL$, since
$$
\rank \bigwedge^t\phi=\binom{\rank\phi}{t}.
$$
If $t=m$, then $Y_t(m,n)=X_t(m,n)$ is just the affine cone over the
Grassmannian. In the case $1<t<mm$ in which we are interested, $Y_t$ is
a proper subset of $X_t$, as we will see soon.

We introduce a numerical invariant that is invariant under the
action of $G=\GL(V)\times\GL(W)$ on $\LL_t$. (It is actually
invariant under the action of $\GL(V)\times \GL(\bigwedge^t W)$.)

\begin{definition}
The \emph{small rank} $\sr(\psi)$ of $\psi\in\LL_t$ is the maximum
of the ranks of the restrictions of $\psi$ to subspaces $\bigwedge^t
U$ of $\bigwedge^t V$ where $U$ ranges over the subspaces of $V$
that have dimension $\le t+1$.
\end{definition}

We now construct elements in $X_t$ for certain pairs of rank and
small rank. Later on we will show that only these numerical
invariants occur in $X_t$ and that they classify the $G$-orbits in
$X_t$. Therefore the diagonal matrices $\dd_{u,u+k-1}$ constructed
in the proof of the following proposition can serve as normal forms.

\begin{proposition}\label{sr-rank}
There exist elements $x\in X_t(m,n)$ for the following combinations
of small rank and rank:
\begin{align*}
\sr(x)&=\rank x=0,\\
\sr(x)&=\rank x=1,\\
\sr(x)&=2,\dots,t+1,\qquad \rank
x=\binom{\sr(x)+k-1}{\sr(x)-1},\quad
 k=1,\dots,m-t.
\end{align*}
\end{proposition}

\begin{proof} Since $m\le n$, we can identify $V$ with a subspace of $W$, and
$\LL_t(V,V)$ with a subspace of $\LL_t(V,W)$. Rank and small rank do
not change if we extend elements from $\LL_t(V,V)$ to $\LL_t(V,W)$
in a trivial way. Therefore we can assume that $m=n$, identify $V$
and $W$, and consider the elements of $L$ as endomorphisms. Let
$e_1,\dots,e_m$ be a basis of $V$.

For $\sr(x)=\rank x=0$ we choose $x=\dd_{0,0}=0$. For $\sr(x)=\rank
x=1$ we choose $x=\dd_{1,1}=\Lambda_t(\phi)$ where $\phi(e_i)=e_i$,
$i=1,\dots,t$ and $\phi(e_i)=0$ for $i>t$.

Let $2\le u\le t+1$, $1\le k\le m-t$, and set $v=t+1-u$. \CHANGE We will now identify an element $\dd_{u,u+k-1}$ in $X_t(m,n)$  with small rank $u$ and rank $\binom{u+k-1}{u-1}$. We consider
the morphism $\alpha:K^*\to L$, where $\alpha(\kappa)$ is the
diagonal matrix with the entries
$$
\alpha(\kappa)_{ii}=\begin{cases}
                        \kappa^{-(t-v)},& 1\le i\le v,\\
                        \kappa^v,& v+1\le i\le v+u+k-1,\\
                        0,& \text{else}.
\end{cases}
$$
Then $\Lambda_t\circ\alpha$ extends to a morphism
$\bar\alpha:K\to\LL_t$ for which $\bar\alpha(0)$ is a diagonal
matrix $\dd_{u,u+k-1}$ with entries $1$ or $0$ on the diagonal.
Clearly $\dd_{u,u+k-1}$ lies in $X_t$.

Furthermore $\dd_{u,u+k-1}$ has  exactly
$$
\binom{u+k-1}{u-1}
$$
entries equal to $1$ on the diagonal and they sit in the positions with indices $[1\dots
v,I\sep 1\dots v,I]$ where $I$ varies over the $(t-v)$-subsets of
$\{v+1,\dots,v+u+k-1\}=\{v+1,\dots,t+k\}$.
It remains to show that $\sr(\dd_{u,u+k-1})=u$. Consider the
subspace $V'=Ke_1+\dots+Ke_{t+1}$. We identify $\bigwedge^t V'$ with
the subspace generated by the basis elements $e_J=\bigwedge_{j\in J}
e_j$ where $J$ is a $t$-subset of $\{1,\dots,t+1\}$. The linear map
$\dd_{u,u+k-1}$ sends $t+1-v=u$ elements of this basis to
themselves, namely those for which $J$ contains $\{1,\dots,v\}$.
Therefore $\sr(\dd_{u,u+k-1})\ge u$.

For the opposite inequality we choose elements $f_1,\dots,f_{t+1}$
in $V$, and represent them in the basis $e_1,\dots,e_m$:
$$
f_i=\sum_{j=1}^m a_{ij}e_j.
$$
Then the restriction of $\dd_{u,u+k-1}$ to $\bigwedge^t V'$,
$V'=Kf_1+\dots+Kf_{t+1}$, is given by a matrix $A'$ whose entries
are $t$-minors of $A=(a_{ij})$. In a row of $A'$ we find the
$t$-minors of $A$ whose row indices leave out a given index
$i=1,\dots,t+1$ and whose column indices correspond to those
$t$-subsets of $\{1,\dots,t+k\}$ that contain $1,\dots,v$. Such a
matrix has rank $\le t+1-v$.

In fact, the rank is maximal when the entries of $A$ are
indeterminates, and then there exist $v$ linearly independent
relations of the $t+1$ rows of $A'$, given by the columns of $A$
with indices $1,\dots,v$  (with appropriate signs), resulting from
Laplace expansion of a $t$-minor with two equal columns (namely the
$j$-th, $j=1,\dots,v$).
\end{proof}

We can already observe that $Y_t$ is a proper subset of $X_t$ if
$1<t<m$. In fact, let $\phi\in L$. If $\rank\phi<t$, then
$\Lambda_t(\phi)=0$, and if $\rank \phi=t$, then
$=\rank\Lambda_t(\phi)=1$. If $\rank\phi\ge t+1$, then
$\sr(\Lambda_t(\phi))=t+1$.

We need some functions which help us to determine small rank.

\begin{lemma}\label{test-sr}
Let $\delta=[1\dots v\sep 1\dots  v]$, $0\le v\le t+1$, and
$\eta=[1\dots t+1\sep1\dots t+1]$ (with $\delta=1$ if $v=0$) and set
$f_v=\delta\eta^{t-v}\in A_t(m,n)=\OO(X_t(m,n))$. Then
$$
f_v(\dd_{u,u+k-1})=0 \quad \iff\quad u<t+1-v.
$$
\end{lemma}

\begin{proof}
First we have to express $f_v$ in the coordinates of $\LL$. We
claim that
$$
f_v=\det\bigl((-1)^{i+j}[1\dots t+1\setminus i\sep 1\dots
t+1\setminus j]_{i,j=v+1,\dots,t+1}\bigr).
$$
Note that this equation generalizes the formula for the
determinant of the adjoint matrix (which it contains for $v=0$).

It is enough to prove the equation over the field of complex
numbers. Both sides of the equation are invariant under the action
of the direct product of the unipotent lower triangular subgroup of
$\GL_m(K)$ and the unipotent upper triangular subgroup of
$\GL_n(K)$. Furthermore they have the same degrees with respect to
all rows and columns of $X$. The space of such forms is
$1$-dimensional (for example, see \cite[(11.11)]{BV}), and so both
sides must differ by a scalar. That it is $1$, follows if we
evaluate both sides on the unit matrix.

Finally we evaluate $f_v$ on the elements $\dd_{u,u+k-1}$.
\end{proof}

\begin{remark}
The results of this section do not depend on characteristic.
\end{remark}

\section{G-stable prime ideals}\label{Atstruct}

The polynomial ring $K[X]\iso\OO(\Hom_K(V,W))$ decomposes into
irreducible $G$-submodules $M_\lambda$ parametrized by the Young
tableaux of shape $\lambda$ as discussed in \cite{DEP} or
\cite[Section 11]{BV}. Each $M_\lambda$ is generated as a $G$-module
by every (standard) bitableau of shape $\lambda$ that is nested on
one side (rows or columns). In particular the bi-initial (or mixed
initial/final) bitableaux belong to $M_\lambda$ and generate it.

The shapes $\lambda$ are non-increasing sequences
$(\lambda_1,\dots,\lambda_u)$ of positive integers such that
$\lambda_i\le\min(m,n)$. We consider some functions on the set of
shapes, namely
$$
\gamma_j(\lambda)=\sum_{i=1}^u \max(\lambda_i-j+1,0),\qquad
j=1,\dots,m,
$$
and
$$
\pi_j(\lambda)=\gamma_j(\lambda)-\frac{\gamma_1(\lambda)}{t}(t-j+1),\qquad
j=1,\dots,t.
$$

Now let $\Pi$ be a product of minors of shape $\lambda$. Then we set
$\gamma_j(\Pi)=\gamma_j(\lambda)$ and $\pi_j(\Pi)=\pi_j(\lambda)$.
The functions $\gamma_j$, introduced in \cite{DEP}, extend to
discrete valuations on $\QF(K[X])$ with non-negative values on
$K[X]$ and the center of $\gamma_j$ is $I_j(X)$ (see \cite{BC2}).
Note that $\gamma_1$ is just the ordinary total degree in $K[X]$.

Note that the $\pi_j$ depend on the value of $t$ under
consideration. The value of $\pi_j$ is an integer for all $j$ if and
only if $\gamma_1(\lambda)$, i.e. the number of boxes of $\lambda$,
is divisible by $t$. In fact, the functions $\pi_j$ are discrete
valuations on the quotient field of the Veronese subalgebra $V_t$ of
the polynomial ring $K[X]$.  See \cite{BC2} for a precise
discussion. There we have shown:  \medskip

\begin{theorem}\label{At}
The subalgebra $A_t$ of $K[X]$ has a basis of standard bitableaux.
One has $A_t=\{x\in V_t:\pi_2(x)\ge 0\}$.
\end{theorem}

\CHANGE The theorem holds for all values of $t,m,n$. For $t=1$ it
holds vacuously since $\pi_2$ is not defined. Note that $V_t$
contains no elements with $\pi_2(x)>0$ if $t=m=\min(m,n)$. If
$1<t<m$, then $V_t$ contains elements of positive value under
$\pi_2$ as, for instance, $\delta^t$ where $\delta$ is a
$(t+1)$-minor of $X$.

In the following we want to work with the weight of a shape (or a
product of minors with that shape). We set
$$
\epsilon_i(\lambda)=\#\{j:\lambda_j=i\}\qquad\text{and}\qquad
\epsilon(\lambda)=(\epsilon_1(\lambda),\dots,\epsilon_m(\lambda)).
$$

The proof of Theorem \ref{At} is based on the following
formula that can be checked by direct computation:

\begin{proposition}\label{pi-formula}
For $u=3,\dots,m$ one has
$\pi_u=(u-1)\pi_2+\sum_{k=1}^{u-2}\epsilon_k$.  In particular,
$\pi_j(x)\ge 0$ for $j\ge 3$ and all $x\in A_t$.
\end{proposition}

As a consequence we obtain

\begin{proposition}
\label{At-dec} The subalgebra $A_t=A_t(m,n)$ is a $G$-submodule of $K[X]$ and
the direct sum of those $M_\lambda$ for which
\begin{itemize}
\item[(i)] $t\mid \gamma_1(\lambda)$ and
\item[(ii)] $\pi_2(\lambda)\ge 0$.
\end{itemize}
\end{proposition}

This follows immediately from Theorem \ref{At} since a product of
minors (especially a nested bitableaux) belongs to $A_t$ if and
only if its degree is divisible by $t$ and its shape $\lambda$
satisfies $\pi_2(\lambda)\ge 0$.

For a $G$-submodule $H$ of $K[X]$ we set
$$
\supp(H)=\{\epsilon(\lambda):M_\lambda\subset H\}.
$$
Clearly $S_t=\supp(A_t)$ is a (normal) submonoid of $\ZZ^m$. It
generates a cone $\RR_+S_t$ in $\RR^m$.

\begin{proposition}
Let $\pp\subset A_t$ be a ($G$-stable) prime ideal. Then
$S_t\setminus\supp(\pp)=S_t\cap F$ where $F$ is a face of
$\RR_+S_t$.
\end{proposition}

\begin{proof}
Consider the subalgebra $B$ of $A_t$ generated by all bi-initial
bitableaux. It is a subalgebra isomorphic to the monoid algebra
$K[S_t]$. Furthermore $\pp\cap B$ is a prime ideal in $B$, and the
ideal $\pp'$ generated by all the monoid elements in $\pp$ (in other
words, the initial bitableaux in $\pp$) is again a prime ideal. But
$\pp'$ is then generated by the monoid elements in $S_t\setminus F$
for some face $F$ of $\RR_+S_t$.
\end{proof}

We denote the face $F$ appearing in the proposition by $\FF(\pp)$.
Now suppose that $\pp$ is $G$-stable. Then $\pp$ is uniquely
determined by $\FF(\pp)$ since $\pp=\Dirsum_{\epsilon(\lambda)\notin
F} M_\lambda$.

We will use the following connection between $G$-stable prime
ideals.

\begin{proposition}
Let $X$ be an affine $K$-variety, and let $G$ be a connected group
acting regularly on $X$. Then the assignment $Gx\mapsto \ol{Gx}$
yields a bijection between the set $\{Gx: x\in X\}$ of orbits and
the set $\{ \ol{Gx} : x\in X\}$ of orbit closures.

If there exist only finitely many orbits or only finitely many
$G$-stable prime ideals in $\mathcal O(X)$, then both these sets are
in bijective correspondence with the set of $G$-stable prime ideals
in $\mathcal O(X)$ via the assignment $Gx\mapsto I(Gx)$.
\end{proposition}

The first assertion follows immediately from the fact that each
orbit closure contains exactly one dense orbit, since each orbit is
open in its closure (see Steinberg \cite{St}). The second assertion
is likewise easily proved.

Let us come back to our variety $X_t$ and its coordinate ring, the
algebras of minors $A_t$. The discussion above shows that $A_t$ has
only finitely many $G$-stable prime ideals.

It is useful to consider the cases $t=1$ and $t=m$ first (we always
assume $m\le n$). If $t=1$, then $A_1=K[X]$ has exactly $m+1$
$G$-stable prime ideals corresponding to the potential ranks of
$m\times n$ matrices, and these are given by the determinantal
ideals $\qq_i=I_i(X)$, $i=1,\dots,m+1$ (with $\qq_{m+1}=0$).

In the case $t=m$ there exist exactly two $G$-stable prime ideals,
namely $\qq_{t+1}=0$, and $\qq_t=I_t(X)\cap A_t$ where the latter is
the irrelevant maximal ideal. This follows immediately from the
transitivity of the action of $G$ on the Grassmannian, but is also
follows from the combinatorial condition on $G$-stable prime ideals
in terms of the monoid ring described a above: if $t=m$, then
$K[S_t]$ is just the polynomial ring in $1$ variable, and it has
only two prime ideals generated by monomials.

These two cases being out of the way, we may assume that $1<t<m$
until we reach Theorem \ref{all-primes}. Under this assumption, the
cone $\RR_+S_t$, introduced above, has exactly $m+1$ facets, namely
$$
F_i=\{(e_1,\dots,e_m)\in \RR_+S_t: e_i\ge 0\},\qquad i=1,\dots,m,
$$
and
$$
F_0=\{e\in \RR_+S_t: \pi_2(e)\ge 0\}.
$$

\begin{proposition}\label{primes1}
The following ideals in $A_t$ are $G$-stable and prime:
\begin{align*}
  &\textup{(1)}\quad \pp_i=\{x\in A_t:\pi_{i+2}(x)> 0\},\quad
  i=0,\dots,t-2;\\
  &\textup{(2)}\quad \qq_j=I_j(X)\cap A_t=\{x\in A_t:\gamma_j(x)>
  0\},\quad  j=t,\dots,m;\\
  &\textup{(3)}\quad \pp_i+\qq_j,\quad i=0,\dots,t-2, j=t+2,\dots,m.
\end{align*}
\end{proposition}

\begin{proof}
The ideals $\pp_i$ and $\qq_j$ are centers of valuations on $A_t$.
Therefore they are prime. Moreover, they are $G$-stable since they
are defined in terms of the $G$-invariant valuations $\gamma_j$.

The best way to show that $\pp_i+\qq_j$ is prime, is to develop the
theory also in the relative situation: one considers the subalgebra
$A_{t,j}$ of $R'=K[X]/I_{j}(X)$ generated by the residue classes of
the $t$-minors. Then $\pi_{i+2}$ (whose definition by shape does not
change in the relative situation) defines a prime ideal in
$A_{t,j}$, and $\pp_i+\qq_j$ is the preimage. (It is not difficult
to transfer \cite{BC2} to the relative version.) An alternative
proof is given below.
\end{proof}

Note that $\qq_t$ is the irrelevant maximal ideal of $A_t$. In
addition to the $t(m-t)+1$ prime ideals listed in Proposition
\ref{primes1} we have the zero ideal, and altogether we have found
$t(m-t)+2$ $G$-stable prime ideals. This is the number of pairs
$(\sr,\rank)$ appearing in Proposition \ref{sr-rank}. For the
following it is useful to set
$$
\pp_{-1}=\qq_{m+1}=0.
$$

\begin{proposition}\label{primes-faces}
For all $i=-1,\dots,t-2$ and $j=t+2,\dots,m+1$ one has
$$
\FF(\pp_i+\qq_j)=F_0\cap\dots\cap F_i\cap F_j\cap\dots\cap F_m.
$$
(where the empty intersection is the full cone). Furthermore
$\FF(\qq_{t+1})=F_0\cap\dots\cap F_{t-1}\cap F_{t+1}\cap\dots\cap
F_m$, and $\FF(\qq_t)=F_0\cap\dots\cap F_m$.
\end{proposition}

\begin{proof}
Since $\qq_t$ is the irrelevant maximal ideal, one has
$\FF(\qq_t)=\{0\}=F_0\cap\dots\cap F_m=\{0\}$.

Now suppose that $\Sigma$ is an bi-initial bitableau whose weight is
not contained in $F_0\cap\dots\cap F_{t-1}\cap F_{t+1}\cap\dots\cap
F_m$. If $\epsilon_i(\Sigma)>0$ for some $i>t$, then clearly
$\Sigma\in\qq_{t+1}$. But if $\epsilon_j(\Sigma)>0$ for some $j<t$,
then $\pi_2(\Sigma)\ge0$ implies that $\epsilon_i(\Sigma)>0$ for
some $i>t$ as well. Since $\FF(\qq_{t+1})$ is properly contained in
$\FF(\qq_t)$, the claim about $\FF(\qq_{t+1})$ follows.

Next we consider $\pp_0$. By definition, $\FF(\pp_0)\subset\FF_0$.
In order to show that $F_j\not\supset\FF(\pp_0)$ for $j>0$ it is
enough to find a bi-initial bitableau $\Sigma\notin\pp$,
equivalently $\pi_2(\Sigma)=0$, with $\epsilon_j(\Sigma)>0$. First
let $j<t$. Then we consider a suitable product $\delta\eta^{t-j}$
where $\delta$ has size $j$ and $\eta$ has size $t+1$. For $j=t$ we
simply take a $t$-minor, and for $j>t$ we consider
$\delta\eta^{j-t}$ where $\eta$ now has size $t-1$.

Similar arguments (together with Proposition \ref{pi-formula})
work in the other cases.
\end{proof}

It follows that there can be no $G$-stable prime ideal strictly
between $\qq_{t-1}$ and $\qq_t$, since there is no face strictly
between $\FF(\qq_{t})$ and $\FF(\qq_{t+1})$.

We want to show that the $G$-stable prime ideals found so far are
the only ones. For this purpose we need the following lemma
\cite[(10.10)]{BV}. In a sense, it describes an anti-straightening
algorithm. It is the basic argument on which Theorem \ref{At} is
based.

\begin{lemma}\label{anti-str}
Let $\delta=[a_1\dots a_u\sep b_1\dots b_u]$ and $\eta=[a_1\dots
a_v\sep b_1\dots b_v]$ with $u<v-1$. Then $\delta\eta$ is a
$K$-linear combination of the products
$$
[a_1\dots a_u a_k\sep b_1\dots b_u b_l][a_1\dots a_v\setminus
a_k\sep b_1\dots b_v\setminus b_l],\qquad k,l=u+1,\dots,v.
$$
\end{lemma}

In proving the converse to Proposition \ref{primes1} we first
characterize the $G$-stable prime ideals different from $\qq_{t+1}$
and $\qq_t$.

\begin{lemma}\label{pimes2}
Let $\pp$ be a $G$-stable prime ideal not containing
$M_{(t+1,t-1)}$. Then $\pp$ is one of the $\pp_i+\qq_j$,
$i=-1,\dots,t-2$, $j=t+2,\dots,m+1$.
\end{lemma}

\begin{proof}
Clearly $\FF(\pp)\not\subset F_t$. In fact, if $\FF(\pp)\subset
F_t$, then all $t$-minors lie in $\pp$, and so $\pp=\qq_t$, the
irrelevant maximal ideal. Moreover, by hypothesis, neither
$\FF(\pp)\subset F_{t-1}$ nor $\FF(\pp)\subset F_{t+1}$. Hence all
the facets containing $\FF(\pp)$ are among
$F_0,\dots,F_{t-2},\allowbreak F_{t+2},\dots,F_m$.

Suppose first that $\FF(\pp)\subset F_i$ for some $i=1,\dots,t-2$.
We have to show that $\FF(\pp)\subset F_{i-1}$ as well, or,
equivalently, that $M_\lambda\subset\pp$ for all $\lambda$ such that
$\epsilon_{i-1}(\lambda)>0$.

By hypothesis none of the row-nested or column-nested bitableaux of
shape $(t+1,\allowbreak t-1)$ is contained in $\pp$, and this will
be very helpful.

Let $\Delta=\delta_1\cdots \delta_w$, $|\delta_1|\le\dots\le
|\delta_w|$ be a bitableau in $A_t$ containing a factor of size
$i-1$, say $|\delta_u|=i-1$. It is certainly enough to show that
$\Delta\in\pp$.

Suppose first that $\Delta$ contains a factor of size $i$. If
$\Delta\in\pp$, we are done. If $\Delta\notin\pp$, the same is true
for all $G$-conjugates of $\Delta$, in particular for those
bitableau produced from $\Delta$ by row permutations. So any product
of such conjugates does not belong to $\Delta$. However a suitable
product of conjugates can be factored into a row-nested bitableau
that has the same shape as $\Delta$ and further factors all of which
also have the same shape. Since the first factor belongs to $\pp$
and the remaining ones are in $A_t$, we obtain a contradiction.

Now suppose there is no factor of size $i$ in $\Delta$. Since
$i-1<t$ and $\pi_2(\Delta)\ge 0$, the product $\Delta$ must contain
a factor, say $\delta_v$, of size $>t$.

We apply anti-straightening to $\delta_u\delta_v$, writing $\Delta$
as a linear combination of products
$$
\Delta'=\delta_1\cdots
\delta_{u-1}\delta'\delta_{u+1}\cdots\delta_{v-1}\delta''
\delta_{v+1}\cdots\delta_{w}
$$
where $|\delta'|=i$, $|\delta''|=|\delta_v|-1$. Now choose a
row-nested bitableau $\zeta\Theta$ of size $(t+1,t-1)$ such that the
rows of the factor $\zeta$ of size $t+1$ contain the rows of
$\delta'$. After multiplication with $\zeta^{t-i}\Theta^{t-i}$
(which is not in $\pp$) we split
$\Delta''=\Delta'\zeta^{t-i}\Theta^{t-i}$ into $\delta'\zeta^{t-i}$
and the product $\Delta'''$ of the remaining factors. Now we have
reached a row-nested bitableau, namely $\delta'\zeta^{t-i}$, with
$\pi_2(\delta'\zeta^{t-i})=0$ and $\gamma_1(\delta'\zeta^{t-i}))$
divisible by $t$. It belongs to $\pp$ since $\FF(\pp)\subset F_i$.
Moreover, the complementary factor $\Delta'''$ belongs to $A_t$
since it has the same value under $\pi_2$ as $\Delta$ and
$\gamma_1(\Delta''')$ is divisible by $t$.

The remaining argument for the case $\FF(\pp)\subset F_j$, $j\ge
t+2$, is almost completely analogous, with $\zeta$ replaced by
$\Theta$. The only exception is that $\Delta$ may have only factors
of size $\ge t$. Then, if $\Delta$ has a factor of size $>t$, one
has $\pi_2(\Delta)>0$. In this case anti-straightening is Laplace
expansion, which reduces the $\pi_2$-value by $1$. But since
$\pi_2(\Delta)>0$, this step is harmless, and the rest of the
argument remains unchanged.
\end{proof}

The remaining case of $\qq_t$ and $\qq_{t+1}$ is handled by the next
lemma.

\begin{lemma}\label{primes3}
Let $\pp\neq \qq_{t}$ be a $G$-stable prime ideal containing
$M_{(t+1,t-1)}$. Then $\pp=\qq_{t+1}$.
\end{lemma}

\begin{proof}
Let $\Delta=\delta_1\cdots\delta_v$ be a product of minors in
$\qq_{t+1}$. We assume that $|\delta_1|\le\dots\le|\delta_w|$.
Moreover, by inserting the empty minor $1$ as an extra factor, we
can assume that $\Delta$ contains a factor of size $\le t-1$.

We want to show that $\Delta\in\pp$, and for an inductive argument
we introduce the following measure:
$$
w(\Delta)=\min\{|\delta_j|-|\delta_i|:|\delta_i|\le t-1,\
|\delta_j|\ge t+1\}.
$$
Note that $|\delta_v|\ge t+1$. Otherwise $\pi_2(\delta)\ge 0$ forces
$|\delta_i|=t$ for all $i$ with $\delta_i\neq 1$, and
$\delta\notin\qq_{t+1}$.

Suppose that $w(\Delta)=2$. Then $\Delta$ contains a factor
$\delta_i\delta_j$ with $\delta_i=t-1$, $\delta_j=t+1$. It may not
be row-nested, but if $\Delta\notin\pp$, then
$\Delta\Delta'\notin\pp$ where $\Delta'$ is a conjugate of $\Delta$
under permutation of the rows. For suitable $\Delta'$ the product
$\Delta\Delta'$ contains a row-nested factor $\eta\zeta$ belonging
to $M_{(t+1,t-1)}$. Since $\pi_2(\eta\zeta)=0$ and
$\gamma_1(\eta\zeta)=2t$, the complementary factor of
$\Delta\Delta'$ is in $A_t$, and we are done.

Now suppose that $w(\Delta)>2$. We choose a pair $\delta_i,\delta_j$
such that $w(\Delta)=|\delta_j|-|\delta_i|$.

Apply anti-straightening to it, and write $\Delta$ as a linear
combination of products $\Theta$ of minors in which $\delta_i$ is
replaced by a minor $\eta$ of size $|\delta_i|+1$ and $\delta_j$ is
replaced by a minor $\zeta$ of size $\delta_j|-1$.

Note that all the products $\Theta$ belong to $A_t$. In fact, only
if $\delta_i=1$, the value under $\pi_2$ drops, and
$\pi_2(\Theta)=\pi_2(\Delta)=-1$ in this case. But then we have
started with $\pi_2(\Delta)>0$.

If $\Theta\notin\pp$, we multiply it by $\Delta$ (if $\Delta\in\pp$,
there was nothing to show). Since not both $\eta$ and $\zeta$ can
have size $t$ (otherwise we had had $w(\Delta)=2$),
$w(\Theta\Delta)<w(\Delta)$, and we are again done.
\end{proof}

In the next theorem the cases $t=1$ and $t=m$ are included again.
However, in these cases we do not define the ideals $\pp_i$ except
$\pp_{-1}=0$.

\begin{theorem}\label{all-primes} There exist exactly $t(m-t)+2$
$G$-stable prime ideals in $A_t(m,n)$, namely
\begin{align*}
  &\textup{(1)}\quad \pp_i+\qq_j,\quad i=-1,\dots,t-2,
  j=t+2,\dots,m+1,\\
  &\textup{(2)}\quad \qq_{t+1}\quad\text{and}\quad \qq_t,
\end{align*}
 where $\pp_{-1}=\qq_{m+1}=0$.
\end{theorem}

\begin{proof}
That the theorem holds in the cases $t=1$ and $t=m$ has been
discussed separately (in these cases $\pp_0=0$).

For $1<t<m$ the preceding results show that only the listed ideals
can be prime and $G$-stable. Therefore there exist at most
$t(m-t)+2$ orbits. That they are prime follows (i) from Proposition
\ref{primes1}, or (ii) from the fact that we must have at least
$t(m-t)+2$ $G$-orbits.
\end{proof}

\begin{corollary}\label{all-pairs}
The pairs of values for small rank and rank listed in Proposition
\ref{sr-rank} are exactly those occurring in $X_t(m,n)$. Each of
them determines a single $G$-orbit.

Moreover, with respect to suitable bases in $V$ and $W$ and the
induced bases of the exterior powers, each element of $X_t(m,n)$ is
given by one of the diagonal matrices $\dd_{u,u+k-1}$ constructed in
the proof of Proposition \ref{sr-rank}.

The $G$-stable prime ideal corresponding to the orbit $\{0\}$ is
$\qq_t$, and $\qq_{t+1}$ corresponds to the orbit of rank $1$
elements in $X_t(m,n)$. The prime ideal corresponding to the orbit
of elements of small rank $u$, $2\le u\le t+1$, and rank
$\binom{u+k-1}{u-1}$, $k=1,\dots,m-t$, is $\pp_{t-u}+\qq_{t+1+k}$
(where again $\pp_{-1}=\qq_{m+1}=0$).
\end{corollary}

\begin{proof}
There are $t(m-t)+2$ such pairs of values, and the number of orbits
is also $t(m-t)+2$.

The $G$-stable prime ideal $\pp$ corresponds to the orbit $Gx$ if
and only if $I(\ol{Gx})=\pp$. It is obvious that the orbit $\{0\}$
and $\qq_t$ correspond to each other. In order to show that
$\qq_{t+1}$ and the rank $1$ matrices correspond to each other, we
can use the element $f_{t-1}$ constructed in Lemma \ref{test-sr}. As
an initial bitableau of shape $(t+1,t-1)$, it belongs to $\qq_{t+1}$
and vanishes only on the rank $1$ element $\dd_{1t}$, Thus only the
orbit of rank $1$ matrices can be contained in $V(\qq_{t+1})$ (in
addition to $\{0\}$).

Using similarly the function $f_v$ which belongs to $\pp_{v}$ for
$v=0,\dots,t-2$ (and counting the $G$-stable prime ideals contained
in $\pp_v$), we see that exactly the orbits given by small rank at
most $t-v$ are contained in $V(\pp_v)$. It follows that only one of
the prime ideals $\pp_{t-u}+\qq_j$, $j=t+2,\dots,m+1$, can
correspond to an orbit with small rank $u$, and now it is enough to
order these orbits by the inclusion of their closures and compare
them to the sequence of prime ideals $\pp_{t-u}+\qq_j$.
\end{proof}

\begin{remark}\label{char-sect2}
(a) Theorem \ref{At}, Lemma \ref{anti-str} and Proposition
\ref{primes1} hold in all non-exceptional characteristics, i.e.\
$\chara K=0$ or $\chara K>\min(t,m-t,n-t)$. The hypothesis that
$\chara K=0$ enters where we have used that a $G$-stable prime ideal
$\pp$ is uniquely determined by $\FF(\pp)$.

(b) The structure of the algebras $A_t$ has been investigated in
\cite{BC2} via a toric deformation. They are Cohen-Macaulay normal
domains in all non-exceptional characteristics, and Gorenstein if
and only if $t=1$, $t=\min(m,n)$, $t=m-1=n-1$, or $1/m+1/n=1/t$.

In arbitrary characteristic, the algebra $\{x\in V_t:\pi_2(x)\ge
0\}$ is ``only'' the normalization of $A_t$; see \cite{B}.

(c) Instead of the toric deformation one can also use the
deformation to the algebra of $U$-invariants where $U$ is a maximal
unipotent subgroup of $G$. For $A_t$ itself this is just the monoid
algebra $K[S_t]$, a normal monoid algebra. By the results of
Grosshans  \cite{Gr} this implies that $A_t$ is normal
Cohen-Macaulay in all non-exceptional characteristics (see also
Bruns and Conca \cite{BC3}). For the residue class rings
$A_t/(\pp_i+\qq_j)$ one obtains the same properties since for them
the algebra of $U$-invariants is generated by the bi-initial
bitableaux in the corresponding face of $\RR_+S_t$. In particular
all orbit closures are normal.

\end{remark}

\section{Analysis of $\Lambda_t$}\label{Lambda-t}

It is now clear that the $G$-orbits in the images $Y_t(V,W)$ are
those of ranks $0$ and $1$ and those of small rank $t+1$. First we
discuss how $\Lambda_t$ acts on the open set of linear maps of rank
$>t$. In order that this set be non-empty we must assume that $t<m$.

\begin{proposition}\label{image1}
Suppose that $t<m$ and let $L_{>t}$ be the open set of linear maps
in $L$ of rank $>t$, and $f,g\in L_{>t}$. Then
$\Lambda_t(f)=\Lambda_t(g)$ if and only if $f=\zeta g$ where $\zeta$
is a $t$-th root of unity. Moreover, if $t>1$, then $\Lambda_t(L_{>t})$ is the complement of
$V(\pp_0)$ in $X_t(m,n)$.
\end{proposition}

\begin{proof}
It follows from Theorem \ref{At} that $(A_t)_x=(V_t)_x$ for all
elements $x$ in $A_t$ with $\pi_2(x)>0$. Therefore the map from the
affine Veronese variety to $X_t(V,W)$ is an isomorphism on the
preimage of $X_t(V,W)\setminus V(\pp_0)$, and thus elements $f,g$ in
$\Lambda_t^{-1}\bigl(X_t(V,W)\setminus V(\pp_0)\bigr)$ go to the
same element in $X_t(V,W)$ if and only if they go to the same
element in the Veronese variety, in other words, if they differ by a
$t$-th root of unity.

Furthermore $X_t(m,n)\setminus V(\pp_0)$ is exactly the union of the
orbits of small rank $>t$. Their union, however, is exactly
$\Lambda_t(L_{>t})$. In particular the fiber over $\Lambda_t(f)$ is
isomorphic to the group of $t$-th roots of unity.

The last statement is clear since among the $G$-stable prime ideals
only $\qq_{t+2},\dots,\qq_{m+1}$ do not contain $\pp_0$, and they
correspond to the orbits of small rank $t+1$.
\end{proof}

\begin{proposition}\label{image2}
Let $L_{t}$ be the set of linear maps in $L$ of rank $t$, and
$f,g\in L_{t}$. Then $\Lambda_t(f)$ and $\Lambda_t(g)$ differ by a
non-zero scalar if and only if $\Ker f=\Ker g$ and $\Im f=\Im g$.
The fiber over $\Lambda_t(f)$ is isomorphic to $\SL_t(K)$.
\end{proposition}

\begin{proof}
Suppose first that $\Ker f=\Ker g$ and $\Im f=\Im g$. Then $f$ and
$g$ both factor through $V/\Ker f$, and can differ only by an
isomorphism $V/\Ker f\to \Im f$. After an identification $V/\Ker
f\iso \Im f$, both $f$ and $g$ can be treated as endomorphisms of
this vector space.  They have the same determinant if and only if
they differ by an element of $\SL(V/\Ker f)\iso \SL_t(K)$.

Now suppose that $\Lambda_t f=\bigwedge^t f$ and $\Lambda_t
g=\bigwedge^t g$ differ only by a nonzero scalar. Choose an element
$x_1\wedge\dots\wedge x_t$ such that $f(x_1)\wedge\dots\wedge
f(x_t)\neq 0$. Set $y_i=g(x_i)$. Then $y_1\wedge\dots\wedge y_t$
belongs to $\Im \bigwedge^t g$ and, in fact, generates it. It
follows that $f(x_1)\wedge\dots\wedge f(x_t)$ and
$y_1\wedge\dots\wedge y_t$ differ only by a scalar. So the subspaces
generated by $y_1\dots,y_t$ and $f(x_1),\dots,f(x_t)$ coincide. By
the dual argument we see that $\Ker f=\Ker g$.
\end{proof}

\begin{proposition}
Suppose that $1<t<m$. Then $\Lambda_t(L_{\le t})$ is the
intersection of $\Im\Lambda_t$ and $V(\pp_0)$. Moreover, the image
of $L_{\le t}$ is closed in $\LL_t$.
\end{proposition}

\begin{proof}
Of the image orbits only those of rank $0$ and rank $1$ can be
contained in $V(\pp_0)$, as seen above. On the other hand,
$\qq_{t+1}$ contains $\pp_0$. The image of $L_{\le t}$ is closed
since it just $V(\qq_{t+1})$.
\end{proof}

\begin{remark}\label{char-sect3}
The results in this section hold in all non-exceptional
characteristics.
\end{remark}

\section{The structure of the orbits}\label{orbits}

Each $\alpha\in V^*$ operates as a derivation on the exterior
algebra $\bigwedge V$ via
$$
x_1\wedge\dots \wedge x_u\hak\alpha=\sum_{i=1}^u (-1)^{i-1}
\alpha(x_i)x_1\wedge\dots\wedge x_{i-1}\wedge x_{i+1}\wedge \dots
\wedge x_u.
$$
This action of $V^*$ on $\bigwedge V$ makes $\bigwedge V$ a right
module over $\bigwedge V^*$ (notation as in Bourbaki \cite[p.\
A.III.162]{Bou:alg}). Furthermore we will consider $\bigwedge W$ as
a left module over itself.

Let
\begin{align*}
D_v(V^*)&=\{\alpha_1\wedge \dots \wedge \alpha_v: \alpha_j\in
V^*\}\setminus \{0\},\\
D_v(W)&=\{y_1\wedge \dots \wedge y_v: y_j\in W\}\setminus \{0\}
\end{align*}
be the subsets of decomposable elements in $\bigwedge^v V^*$ and
$\bigwedge^v W$. Since $D_v(W)$ is the quasiprojective variety of
the nonzero elements in the affine cone over $\Grass_v(W)$. One has
$\dim D_v(W)=v(n-v)+1$, and similarly $\dim D_v(V^*)=v(m-v)+1$.

We define maps
\begin{align*}
&\Theta_{\alpha,y}:\LL_{t-v}(V,W)\to \LL_t(V,W),\qquad
\alpha\in D_v(V^*),\ y\in D_v(W),\\
&\bigl(\Theta_{\alpha,y}(f)\bigr)(x)=y\wedge f(x\hak\alpha).
\end{align*}

Let
$$
V_\alpha=\bigcap_{i=1}^v \Ker\alpha_i\quad\text{and}\quad
W_y=W/\sum_{i=1}^vKy_i.
$$
Note that $(\bigwedge^t V)\hak\alpha=\bigwedge^{t-v}V_\alpha$ and
that $y\wedge z=0$ for all $z$ in the kernel of the natural map
$\bigwedge^{t-v}W \to \bigwedge^{t-v}W_y$. Hence
$\Theta_{\alpha,y}(f)$ only depends on the map
$f'\in\LL_{t-v}(V_\alpha,W_y)$ induced by $f$ (via restriction to
$\bigwedge^{t-v}V_\alpha$ and the composition with the projection
$\bigwedge^{t-v}W \to \bigwedge^{t-v}W_y$). Therefore we will
consider $\Theta_{\alpha,y}$ to be defined on
$\LL_{t-v}(V_\alpha,W_y)$.

\begin{proposition}\label{retracts2}
$$
\Theta_{\alpha,y}\bigl(X_{t-v}(V_\alpha,W_y)\bigr)\subset X_t(V,W)
$$
is a retract of $X_t(V,W)$.
\end{proposition}

\begin{proof}
First we have to show that indeed
$\Theta_{\alpha,y}\bigl(X_{t-v}(V_\alpha,W_y)\bigr)\subset
X_t(V,W)$. Consider an element $\phi\in X_{t-v}(V_\alpha,W_y)$. As
we have seen in Corollary \ref{all-pairs}, with respect to a
suitable bases $e_1,\dots, e_{m-v}$ of $V_\alpha$ and
$f_1,\dots,f_{n-v}$ of $W_y$ and the induced bases on the exterior
powers, the matrix of $\phi$ has the form $\dd_{u,u+k-1}$.

Now choose bases $e_1',\dots,e_m'$ of $V$ and $f_1',\dots,f_n'$ of
$W$ such that $e_i'=e_{i-v}$ for $i=v+1,\dots,m$ and $f_j'=f_{j-v}$
for $j=v+1,\dots,n$. Then $\Theta(\dd_{u,u+k-1})=\dd_{u,u+k-1}$
where the matrix on the right is formed in $\LL_t(V,W)$ with respect
to the bases of $\bigwedge^t V$ and $\bigwedge^t W$ induced by
$e_1',\dots,e_m'$ and $f_1',\dots,f_n'$.

It is not hard to define a map
$\Psi:\LL_t(V,W)\to\LL_{t-v}(V_\alpha,W_y)$ such that
$\Psi\circ\Theta_{\alpha,y}$ is the identity on
$\LL_{t-v}(V_\alpha,W_y)$. We extend $\alpha_1,\dots\alpha_v$ to a
basis $\alpha_1,\dots,\alpha_m$ of $V^*$ and $y_1\dots,y_v$ to a
basis of $W$. Then we take the dual bases $(\alpha_i^*)$ of $V$ and
$(y_j^*)$ of $W^*$. Now let
$$
(\Xi(\phi))(x)=(y_1^*\wedge\dots\wedge
y_v^*)\bigl(\phi(\alpha_1^*\wedge\dots\wedge \alpha_v^*\wedge
x)\bigr).
$$
It remains to show that $\Xi(X_t(V,W))\subset
X_{t-v}(V_\alpha,W_y)$. This follows from the next proposition.
\end{proof}

We want to write out the maps defined in terms of coordinates. To
this end we let (after the choice of bases as above) denote the
coordinates on $\LL_t(m,n)$ by $E_{I,J}$ where $I\subset\{1,\dots,m\}$ and
$J\subset\{1,\dots,n\}$, $\#I=\#J=t$. By $E'$ we denote the corresponding
coordinates of $\LL_{t-v}(m-v,n-v)$. Then $\Theta=\Theta_{\alpha,y}$
is given by the substitution
$$
\theta: E_{I,J}\mapsto\begin{cases}
                    E'_{I^-,J^-},& \{1,\dots,v\}\subset I,J,\\
                    0,&\text{else}.
\end{cases}
$$
Here we set $I^-=\{i_{v+1}-v,\dots,i_t-v\}$ and $J^-$ is defined
accordingly. As we have seen already, this substitution induces a
morphism of coordinate rings $A_t(m,n)\to A_{t-v}(m-v,\allowbreak
n-v)$ sending $[I\sep J]_X$ to $[I^-\sep J^-]_{X'}$ if
$\{1,\dots,v\}\subset I,J$, and to $0$ otherwise.

In order to complete the proof of Proposition \ref{retracts2} we
have to show that the substitution corresponding to $\Psi$, namely
$$
\psi:E'_{L,M}\mapsto E_{L^+,M^+}, \qquad
L^+=\{1,\dots,v\}\cup\{l_1+v,\dots,l_{t-v}+v\}
$$
induces a morphism $A_{t-v}(m-t,n-t)\to A_t(m,n)$.

To this end we consider the polynomial rings $K[X]$ and $K[X']$ and
the substitution
$$
\zeta: X_{ij}\mapsto\begin{cases}
                    1,&i=j\le v,\\
                    X'_{i-v,j-v},& i,j>v,\\
                    0,&\text{else}.
\end{cases}
$$
Let $B$ be the subalgebra of $A_t(m,n)$ generated by the minors
$[L^+\sep M^+]_X$. Evidently $\zeta$ induces a surjective morphism
$B\mapsto A_{t-v}(m-v,n-v)$. We claim that it is an isomorphism.
Then the inverse is exactly the morphism induced by $\psi$, and we
have proved Proposition \ref{retracts2}.

Note that $B$ has a basis $S$ of standard bitableaux all of whose
factors contain $\{1,\dots,v\}$ in their row as well in their column
part. This follows immediately from the straightening algorithm,
since straightening preserves the content of a bitableau and each
generator of $B$ contains $\{1,\dots,v\}$ in its row and column
part. In particular, the standard bitableaux in a representation of
an element of normalized degree $k$ all have exactly $k$ factors.
Clearly $\gamma$ maps a standard bitableau in $S$ of degree $k$ to a
standard bitableau in $A_{t-v}(m,n)$ of the same degree $k$, and in
fact induces a bijection between the degree $k$ elements in $S$ and
the degree $k$ standard bitableau in $A_{t-v}(m-v,n-v)$. To sum up:

\begin{proposition}\label{retracts}
The substitutions $\theta$ and $\psi$ induce a retraction of
$A_t(m,n)$ onto $A_{t-v}(m-v,\allowbreak n-v)$.

More generally, $A_{t'}(m',n')$ is a retract of $A_t(m,n)$ if
$m-m',n-n'\ge t-t'$.
\end{proposition}

The last statement is very easy to see. Set $m''=m'-(t-t')$ and
$n''=n-(t-t')$. Then $A_t(m'',n'')$ is evidently a retract of
$A_t(m,n)$, and then we can apply the first part.

Next we want to investigate the effect of varying $\alpha$ and $y$
in Proposition \ref{retracts2}. To this end we define a suitable
bundle over $D_v(V^*)\times D_v(W)$. First we consider the trivial
bundles with fiber $\bigwedge^{t-v} V$ over $D_t(V^*)$ and fiber
$\bigwedge^{t-v} W$ over $D_t(W)$. Then we take the subbundle of the
first whose fiber over $\alpha_1\wedge\dots\wedge \alpha_v$ is
$V_\alpha$, and the quotient bundle of the second whose fiber over
$y_1\wedge\dots\wedge y_v$ is $W_y$. Finally we get the bundle
$\BB_{v,t-v}$ over $D_v(V^*)\times D_t(W)$ whose fiber over
$(\alpha,y)$ is $\LL_{t-v}(V_\alpha,W_y)$.

The bundles we have defined are locally trivial and quasi-projective
varieties.

\begin{theorem}\label{orb-struct}
Let $u=1$, $k=0$ or $u\in\{2,\dots,t+1\}$ and $k\in\{1,\dots,m-t\}$.
Set $v=t-u+1$ and let $Z$ be the subvariety of $\BB_{v,t-v}$ whose
fiber over $(\alpha,y)$ is the $\GL(V_\alpha)\times \GL(W_y)$-orbit
$O$ of $Y_{u-1}(V_\alpha,W_y)$ corresponding to small rank $u$ and
rank $\binom{u-1+k}{u-1}$. Then its image under
$$
\Theta:\BB_{v,t-v}\mapsto X_t(V,W),\qquad
\Theta(\alpha,y,f)=\Theta_{\alpha,y}(f),
$$
is the $G$-orbit $O'$ in $X_t(V,W)$ with small rank $u$ and rank
$\binom{u-1+k}{u-1}$.

Moreover, $\Theta(\alpha',y',f')=\Theta(\alpha,y,f)$ for $f,f'\in
Y_{t-v}(V_\alpha,W_y)$if and only if there exist
$\lambda,\mu,\rho\in K^*$ such that $\alpha'=\lambda\alpha$, $y'=\mu
y$, $f'=\rho f$ and $\lambda\mu\rho=1$. (For $u=1$, $k=0$ one
necessarily has $\rho=1$.)
\end{theorem}

\begin{proof}
By the definition of $\Theta_{\alpha,y}$, the image of $\Theta$
really is the image of the map $D_v(V^*)\times D_v(W)\times
\LL_{t-v}(V,W)$ that sends $(\alpha,y,f)$ to $x\mapsto y\wedge
f(x\hak\alpha)$. This map is $G$-equivariant, and therefore maps
orbits to orbits. Thus the first statement holds.

We will now show that we can recover $y$ up to a nonzero scalar from
the image of $\Theta(\alpha,y,f)$ if $f\in Y_{t-v}(V_\alpha,W_y)$.
We choose $g\in\Hom_K(V_\alpha,W_y)$ such that $f=\bigwedge^{t-v}g$.
We extend $g$ to $V$ and lift the extension to a linear map, also
called $g$, from $V$ to $W$. Then, as discussed above,
$\Theta_{\alpha,y}(f)$ is given by $x\mapsto y\wedge
f(x\hak\alpha)$. The image of $f\circ\alpha$ is $f(\bigwedge^{t-v}
V_\alpha)$, since $\bigwedge^{t-v} V_\alpha$ is the image of
$\alpha$. Since $f=\bigwedge^{t-v} g$ we obtain finally that the
image of $\Theta_{\alpha,y}(f)$ is of the form
$y\wedge\bigwedge^{t-v} W'$. Since $f$ has small rank $t-v+1$, the
image of $W'$ in $W_y$ has dimension at least $t-v+1$. But then it
follows that the annihilator of $y\wedge\bigwedge^{t-v} W'$ in $W$
with respect to exterior multiplication is the subspace generated by
$y_1,\dots,y_v$. This subspace determines $y=y_1\wedge\dots\wedge
y_v$ up to a nonzero scalar. It follows that $y'=\mu y$ if
$\Theta(\alpha',y',f')=\Theta(\alpha,y,f)$. By dual reasons we
obtain $\alpha'=\lambda \alpha$.

We can replace $y'$ by $y$ and $\alpha'$ by $\alpha$, and,
consequently, $f'$ by $\lambda^{-1}\mu^{-1}f'$. Then
$f=\lambda^{-1}\mu^{-1}f'$ by Proposition \ref{retracts2}.

Note that $Y_0(V_\alpha,W_y)$ contains only one element, namely the
multiplication by $1$ on $K$. This explains why $\rho=1$ in the case
$u=1$, $k=0$.
\end{proof}

\begin{corollary}\label{dimension}
The orbit $O'$ of $rank 1$ in $X_t(m,n)$ has dimension
$(m-t)t+(n-t)t+1$. The orbit $O'$ of small rank $u\ge 2$ and rank
$\binom{u-1+k}{u-1}$, $k=1,\dots,m-t$, has dimension
$$
mn-(t-u+1)^2-\bigl(m-(t+k)\bigr)\bigl(n-(t+k)\bigr).
$$
\end{corollary}

\begin{proof}
By the theorem,
$$
\dim O'=\dim D_v(V^*)+\dim D_v(W)+\dim O -\epsilon
$$
where $v=t+1-u$ and $\epsilon=2$ for $u>1$ and $\epsilon=1$ for
$u=1$. (This discrepancy explains why the formula for $u\ge 2$, when
applied in the case $u=1$, gives a result that is too small by $1$.)

One has $\dim D_v(V^*)=v(m-v)+1$ and $\dim D_v(W)=v(n-v)+1$. It
remains to find the dimension of $O'$. As we have seen above, $O'$
is the quotient of the quasiprojective variety $L_{t-v+k}$ of rank
$u-1+k$ linear maps in $L=\Hom_K(V_\alpha,W_y)\iso
\Hom_K(K^{m-v},K^{n-v})$ modulo the action of the $(u-1)$-th roots
of unity. This is a finite group and so $\dim O'=\dim L_{u-1+k}$.
But $\dim L_{u-1+k}=\dim L_{\le u-1+k}$, and the standard formula
for the dimension of the determinantal variety $L_{\le u-1+k}$
yields
\begin{equation*}
\dim O'=(m-v)(n-v)- ((m-v)-(u-1+k))((n-v)-(u-1+k)).\qedhere
\end{equation*}
\end{proof}

There is a natural isomorphism $\LL_i(V,W)\iso \bigwedge^i
V^*\tensor \bigwedge^i W$, $i\in\NN$ where the action of
$\bigwedge^i V^*$ is the right multiplication on $\bigwedge^i V$
given above. With respect to this isomorphism, the assignment
$\phi\mapsto \bigwedge^i \phi$ is just the $i$th power map in the
algebra $\bigwedge V^*\tensor \bigwedge W$ or its subalgebra
$\Dirsum_t \bigwedge^t V^*\tensor \bigwedge^t W$. The  map
$\Theta_{\alpha,y}$ is the left multiplication by $\alpha\tensor y$:
\begin{align*}
\bigl(\Theta_{\alpha,y}(\beta\tensor
z)\bigr)(x)&=y\wedge(\beta\tensor z)(x\hak\alpha)\\
&=y\wedge((x\hak\alpha\hak\beta)z)\\
&=(x\hak(\alpha\wedge\beta))y\wedge z,
\end{align*}
and we get the same result by applying
$(\alpha\wedge\beta)\tensor(y\wedge z)$ to $x$. (Note that
$x\hak\alpha\hak\beta\in K$.)

If we give up some precision in the description of the orbits, we
get a very smooth result for the Zariski closure of the union
$$
X(V,W)=\bigcup_{t\ge 0} X_t(V,W)
$$
in the algebra $\bigwedge V^*\tensor \bigwedge W$. Analogously we
define $Y(V,W)$ as the union of the images $Y_t(V,W)$. The following
theorem has already been stated in the introduction.

\begin{theorem}
$X(V,W)$ is the closure of $Y(V,W)$ under the iterated operation of
$V^*\times W$ by multiplication on $\bigwedge V^*\tensor \bigwedge
W$.

Moreover, for $x\in X_t(V,W)$ with $\rank x>1$ there exists a
presentation
\begin{multline*}
x=\bigl((\alpha_1\wedge\dots\wedge\alpha_{t-s})\tensor
(y_1\wedge\dots\wedge y_{t-s})\bigr)\cdot x',\\
x'\in Y_s(V,W),\ \alpha_1,\dots,\alpha_{t-s}\in V^*,
y_1,\dots,y_{t-s}\in W,
\end{multline*}
if and only if $s=\sr(x)-1$. (If $\rank x=1$, then one can choose
$s=\sr(x)-1=0$.)
\end{theorem}

\begin{proof}
Only the statement on $\sr(x)-1$ has not yet been completely
justified. First let $s=\sr(x)-1$. To guarantee the existence of a
representation as claimed in the theorem, it is enough to check that
the normal forms $\dd_{u,u+k-1}$ of small rank $u$ indeed ``come
from'' $Y_{u-1}(V,W)$, and this has been seen in the proof of
Proposition \ref{retracts2}.

Suppose now that $\rank x>1$ and that $x$ has a representation as in
the theorem. Then $x=\Theta_{\alpha,y}(x'')$ where $x'\in
\LL_{s}(V_\alpha, W_y)$ is induced by $x'$. Evidently $x''\in
Y_{s}(V_\alpha,W_y)$. Therefore $\sr(x'')$ can only attain the
values $0,1,s+1$. But if $\sr(x'')\le 1$, then $\rank x''\le 1$ as
well, and evidently $\rank x\le 1$. So only $\sr(x'')=s+1$ is
possible, and since $\sr(x)=\sr(x'')$, we conclude that
$s=\sr(x)-1$. (The equation $\sr(x)=\sr(x'')$ can again be checked
on the normal forms.)
\end{proof}

\begin{remark}
The hypothesis that $\chara K=0$ enters the results of this section
only through Corollary \ref{all-pairs}. Again, non-exceptional
characteristic would be enough if Corollary \ref{all-pairs} could be
generalized.
\end{remark}

\section{The singular locus}\label{Sing}

The goal of this section is to identify the singular locus
$\Sing(X_t)$ of $X_t$. This will be achieved by describing some
localizations of $A_t$. We need some preliminary results. The first
is Binet's formula:

\begin{lemma}\label{Binet}
Let $A$ and $B$ be matrices of size $m\times n$ and $n\times p$. Set
$C=AB$. Let $a=a_1,\dots,a_t$ and $b=b_1,\dots,b_t$ with $1\leq
a_i\leq m$ and $1\leq b_i\leq m$. Then
$$[a\sep b]_C=\sum [a\sep c]_A\ [c\sep b]_B$$
where the sum is extended to all the increasing subsequences $c$ of
length $t$ of $\{1,\dots,n\}$.
\end{lemma}

The second is a way of producing new relations among minors of a
given size by starting from known relations and applying the group
operation. A \emph{Pl\"ucker relation} on $t$-minors is a quadratic
relation with integral coefficients among the maximal minors of the
generic $t\times n$ matrix. It can be written in the form:
\begin{equation}
\sum \lambda_i [\alpha_i]\ [\beta_i] =0 \label{general}
\end{equation}
where $\lambda_i \in \ZZ$, and $\alpha_i$, $\beta_i$ are increasing
sequences on length $t$ bounded below by $1$ and above by $n$ and
$[\alpha_i]$ and $[\beta_i]$ are the corresponding $t$-minors. We
say that the Pl\"ucker relation (\ref{general}) is
\emph{homogeneous} if the multi-set $\alpha_i\union \beta_i$ does
not depend on $i$; in that case we say that (\ref{general}) is
homogeneous of \emph{degree} $\alpha_i\union \beta_i$. The typical
homogeneous Pl\"ucker relation arises in the following way: Let
$a=a_1,\dots,a_{t-1}$ and $b=b_1,\dots, b_{t+1}$ column indices;
then the Pl\"ucker relation associated with $a$ and $b$ is:
\begin{equation}
\sum_{j=1}^{t+1} (-1)^{j+1} [a_1\dots a_{t-1}, b_j]\ [b_1\dots
b_{j-1}, b_{j+1} \dots b_{t+1}]=0 \label{typical}
\end{equation}
For example, with $t=2$ and $a=1$ and $b=2,3,4$ we get
\begin{equation}
[12][34]-[13][24]+[14][23]=0 \label{plu2}
\end{equation}
This is essentially the only Pl\"ucker relation on $2$-minors. It is
homogeneous of degree $\{1,2,3,4\}$.

\begin{lemma}\label{genplu1}
Consider a homogeneous Pl\"ucker relation (\ref{general}) on
$t$-minors of degree $v$. For every matrix $X$ of size $m\times n$,
for every matrix $A$ of size $t\times m$ and for every multi-set $u$
of elements in $\{1,\dots,m\}$ and of cardinality $2t$ one has:
$$
\sum_c [c]_A\ [c']_A \sum_i \lambda_i [c\sep \alpha_i]_X\ [c'\sep
\beta_i]_X=0
$$
where the sum $\sum_c$ is extended over all the increasing sequences
$c\subset u$ of cardinality $t$, and $c'=u\setminus c$.
\end{lemma}

\begin{proof} We may assume that $X=(X_{ij})$
is a matrix of variables. We give $X_{ij}$ the multidegree
$(e_i,f_j)$ in $\ZZ^{m}\times \ZZ^{n}$ where $e_i$ and $f_j$ are the
corresponding canonical bases of $\ZZ^{m}$ and $\ZZ^n$. The
Pl\"ucker relation, specialized to the matrix $AX$ gives a relation
among the maximal minors of $AX$. Using Binet's formula we get
$$
\sum_{c,b} [c]_A\ [b]_A \sum_i [c|\alpha_i]_X\ [b| \beta_i]_X=0
$$
where the sum $\sum_{c,b}$ is extended over all the ordered pairs
$(c,b)$ of increasing sequences of length $t$ in $\{1,\dots, m\}$.
So all the multihomogeneous components of the above polynomial in
the $X_{ij}$ vanishes. The vanishing of the component of multidegree
$(r_1,\dots,r_m,c_1,\dots,c_n)$ where $r_k=|\{ i \in u : i=k\}|$ and
$c_k=|\{ i \in v : i=k\}|$ gives the desired expression.
\end{proof}

Particular choices of the matrix, of the Pl\"ucker relation
(\ref{general}), of $A$ and of $u$ result in particular kinds of
relations:

\begin{lemma}
\label{genplu2} For every $0\leq s\leq t$ set $w=\{s+1, s+2,\dots,
2t-s\}$ and consider the set $B$ of the subsequences $b$ of $w$ of
cardinality $t-s$ and such that for all $i=s+1,\dots,s+t$ one has
$|b\sect \{i, i+t-s\}|=1$. Then the relation
$$
\sum_{\bf b\in B} (-1)^b \sum_{j=t}^{2t} (-1)^{j+1} [1\dots s, b
\sep 1\dots t-1 j ]\ [1\dots s, b' \sep t\dots \widehat{j} \dots 2t
] =0
$$
holds in every matrix. Here $b'=w\setminus b$ and $(-1)^b=\pm 1$.
\end{lemma}

\begin{proof} Let $A$ be the $t\times (2t-s)$ matrix with block decomposition
$$
\begin{pmatrix} I_s & 0 & 0 \\ 0 & I_{t-s} & I_{t-s}
\end{pmatrix}
$$
where $I_k$ denotes the $k\times k$ identity matrix. Set
$u=\{1,1,2,2,\dots,s,s,s+1, s+2,\dots, 2t-s\}$. Then apply Lemma
\ref{genplu1} to the Pl\"ucker relation (\ref{typical}) with $A$ and
$u$ as just defined. To see that one gets the claimed relation note
that the non-zero maximal minors of $A$ are all $1$ or $-1$, that
they arise exactly by taking column indices of the form $1\dots s,
b$ with $b\in B$, and that $b\in B$ if and only if $b'\in B$.
\end{proof}

Indeed simple considerations show that the sign $(-1)^b$ in the
formula above is $1$ if $t-s$ is odd or if the cardinality of $\{ i
\in b : i\leq t\}$ is even and $-1$ otherwise.

For instance, by applying \ref{genplu2} with $s=0$ and $t=2$, we
have $B=\bigl\{ \{1,2\}, \{1,4\}, \{2,3\},\allowbreak \{3,4\}
\bigr\}$ and get the $12$-term relation
\begin{equation}
\begin{array}{l}
+[12|14][34|23]-[14|14][23|23]-[23|14][14|23]+[34|14][12|23] \\[6pt]
 -[12|13][34|24]+[14|13][23|24]+[23|13][14|24]-[34|13][12|24] \\[6pt]
 +\underline{[12|12][34|34]}-[14|12][23|34]-[23|12][14|34]+[34|12][12|34]\quad=0.
\end{array}
\label{12terms2min}
\end{equation}
valid in every $4\times 4$ matrix. It will become apparent below why
one of the terms has been underlined.

The retract property given in Proposition \ref{retracts} implies
that appending a new row index and a new column index to a known
relation among $t$-minors yields a relation for $(t+1)$-minors. For
instance, taking relation \eqref{12terms2min}, appending a new row
and column index, say $0$, and then shifting the indices by $1$ we
get the relation
\begin{equation}
+[123|125][145|134]-\dots +\underline{[123|123][145|145]}-\dots
+[145|123][123|145]=0 \label{12terms3min}
\end{equation}
valid in every $5\times 5$ matrix.

So far we have seen relations of degree $2$ which are, modulo the
group operation, Pl\"ucker relations. The following are important
relations of degree $3$ among $2$-minors of any $4\times 4$ matrix:
\begin{subequations}
\begin{align}
\bl  \{12\}, \{13\}, \{24\} \bs  \{12\}, \{13\}, \{24\} \br&= \bl
\{12\},
\{14\}, \{23\} \bs  \{12\}, \{14\}, \{23\} \br \label{6a}\\[6pt]
\bl \{12\},\{13\}, \{23\} \bs  \{12\}, \{13\}, \{24\} \br&= \bl
\{12\},
\{13\}, \{23\} \bs  \{12\}, \{14\}, \{23\} \br \label{6b}\\[6pt]
\bl \{12\}, \{13\}, \{14\} \bs  \{12\}, \{13\}, \{24\} \br&=- \bl
\{12\},
\{13\}, \{14\} \bs \{12\}, \{14\}, \{23\} \br\label{6c} \\[6pt]
\bl \{12\}, \{13\}, \{14\} \bs \{12\}, \{13\}, \{23\} \br&=0 \label{6d}\\[6pt]
[13\sep 24]\cdot\bl\{12\},\{13\} \bs \{12\}, \{13\}\br&=\notag\\[6pt]
[13\sep
23]\cdot & \bl\{12\}, \{13\} \bs \{12\}, \{14\}\br+ [13\sep12]\cdot G\label{6e}\\[6pt]
 \mbox{ with }  G=-[12\sep23 ] [13\sep 14 ]+[12\sep 24 ]&[13\sep 13 ] -[12\sep34 ]
 [13\sep 12]\notag
\end{align}
\end{subequations}
Relations \eqref{6a}--\eqref{6d} are described in terms of
$3$-minors of $2$-minors. For instance, \eqref{6d} says that
$$\det  \begin{pmatrix}
 [12| 12]  &  [12|13] & [12|23]\\
 [13| 12]  &  [13|13] & [13|23]\\
 [14| 12]  &  [14|13] & [14|23]
 \end{pmatrix} =0,
$$
and this should suffice to explain our notation. The relations
\eqref{6a}--\eqref{6e} can be checked directly by expansion (it is a
good idea to use a computer algebra system for this task). Note
however that \eqref{6b} is obtained from \eqref{6a} by replacing the
row index $4$ with $3$, \eqref{6d} is obtained from \eqref{6c} by
replacing the column index $4$ with $3$ and dividing by $2$. So it
is enough to check \eqref{6a}, \eqref{6c} and \eqref{6e}. With some
more effort one can check that \eqref{6a}, \eqref{6b} and \eqref{6c}
arise by applying \eqref{6d} to a matrix of the form $AXB$ with $A$
and $B$ scalar matrices and then selecting homogeneous components;
in other words; they arise from \eqref{6d} by the operation of $G$.
Also note that \eqref{6d} follows immediately from the obvious fact
that the vectors $x_1\wedge x_2,x_1\wedge x_3, x_1\wedge x_4$ are
linearly dependent if $x_1,\dots,x_4$ are so. Also \eqref{6e}
results from \CHANGE   \eqref{6d}, Pl\"ucker relations and their $G$-conjugates  since all relations of the $2$-minors
 arise in this way (see \cite{B}).

Returning to our algebra of minors, we want to find localizations of
$A_t$ which are regular rings. Precisely, we will describe a subset
$\Phi_0$ of $t$-minors such that $\Phi_0$ has cardinality $mn$ and
the $K$-subalgebra $K[\Phi_0]$ of $A_t$ generated by the elements of
$\Phi_0$ coincide with $A_t$ after the inversion of a suitable
element $F$ of $A_t$. Since we need a $(t+1)$-minor of the matrix
$X$, we must assume that $t<m$.

Let us describe $\Phi_0$, $F$ and also two auxiliary sets $\Phi_1$
and $\Phi_2$. Let $\delta_i$ be the minor of the first $i$ rows and
columns of $X$. Then set
$$
F=\delta_{t-1}\delta_t\delta_{t+1}.
$$

By definition, $\Phi_0$, $\Phi_1$ and $\Phi_2$ are the sets of the
$t$-minors $[a_1\dots a_t | b_1\dots  b_t]$ with $1\leq
a_1<\dots<a_t\leq m$ and $1\leq b_1<\dots<b_t\leq n$ defined as
follows. In each of the three cases we require that \emph{at least}
one of the conditions is satisfied:
$$
\begin{aligned}
\Phi_0:\ \text{(i)}&\  a_t\leq t+1 \text{ and } b_t\leq t+1\\
\text{(ii)}&\ a_{t-1}=t-1 \text{ and } b_{t-1}=t-1\quad\\
\text{(iii)}&\ a_{t}=t \text{ and } b_{t-1}\leq t\\
\text{(iv)}&\ a_{t-1}\leq t \text{ and } b_{t}=t\\[6pt]
\Phi_2:\ \text{(i)}&\  a_{t-1}\leq t \text{ and } b_{t-1}\leq t\\
\text{(ii)}&\ a_{t}= t\\
\text{(iii)}&\  b_{t}=t
\end{aligned}
\
\begin{aligned}
\Phi_1:\ \text{(i)}&\  a_t\leq t+1 \text{ and } b_t\leq t+1\\
\text{(ii)}&\ a_{t-1}=t-1 \text{ and } b_{t-1}=t-1\\
\text{(iii)}&\  a_{t}=t\\
\text{(iv)}&\  b_{t}=t\\[6pt]
\vphantom{\Phi_2}\\
\vphantom{a_{t}= t}\\
\vphantom{b_{t}=t}
\end{aligned}
$$
By definition, $\Phi_0\subset \Phi_1 \subset \Phi_2$. Let
$K[\Phi_i]$ denote the subalgebra of $A_t$ generated by the elements
of $\Phi_i$. Note that $F\in K[\Phi_0]$ since
$$\delta_{t-1}\delta_{t+1}=\det\begin{pmatrix}
[1\dots t | 1\dots t]  & [1\dots t |1\dots t-1,t+1] \cr \cr [1\dots
t-1,t+1 | 1\dots t] & [1\dots t-1,t+1 | 1\dots t-1,t+1]
\end{pmatrix}.
$$

We have:

\begin{theorem}
\label{localize} Suppose that $t<m$. Then
\begin{itemize}
\item[(1)] the cardinality of $\Phi_0$ is $mn$;
\item[(2)] $K[\Phi_0][F^{-1}]=A_t[F^{-1}]$;
\item[(3)] $A_t[F^{-1}]$ is a regular ring.
\end{itemize}
\end{theorem}

\begin{proof}
(1) is a simple count, taking care of the overlaps in the conditions
defining $\Phi_0$. (3) follows from (1) and (2). Indeed, (2) implies
that $K[\Phi_0]$ and $A_t$, having the same field of fractions, have
the same dimension $mn$. Combined with (1), we obtain that
$K[\Phi_0]$ is a polynomial ring, and so (3) follows. Hence the
crucial statement is (2). We prove it in three steps: first we show
that $K[\Phi_2][F^{-1}]=A_t[F^{-1}]$, then that
$K[\Phi_1][F^{-1}]=K[\Phi_2][F^{-1}]$ and finally that
$K[\Phi_0][F^{-1}]=K[\Phi_1][F^{-1}]$ . This is done in the
following lemma.
\end{proof}

\begin{lemma}
\label{criticalClaim} Suppose that $t<m$. Then
\begin{itemize}
\item[(1)] for every $t$-minor $M$ of $X$ there exists an integer
$k$ such that $\delta_t^kM\in K[\Phi_2]$;
\item[(2)] for every $t$-minor $M$ in $\Phi_2$ there exists an integer
$k$ such that $(\delta_{t-1}\delta_{t+1})^kM\in K[\Phi_1]$;
\item[(3)] for every $t$-minor $M$ in $\Phi_1$ there exists an integer
$k$ such that $\delta_t^kM\in K[\Phi_0]$.
\end{itemize}
\end{lemma}

\begin{proof}
(1) Let $M$ be a $t$-minor. Let $e(M)=(s,v)$ where $s$ is the number
of row column indices of $M$ which which are $\leq t$ and $v$ is the
number of column indices of $M$ which are $\leq t$. If $s=t$ or
$v=t$ or $v=s=t-1$, then $M\in \Phi_2$. We now argue by decreasing
induction on $t$ and on $(s,v)$:

Case (a): if $s>0$ and $v>0$ then we may assume that $M$ involves
the first row and column. Then, by the principle of retraction, we
can reduce the statement to the case of $(t-1)$-minors and are done
by induction.

Case (b): if $s$ or $v$ is $0$, we may assume $v=0$, transposing if
necessary. Then, up to a renaming of the indices, we may further
assume that
$$
M=[1\dots s,t+1\dots  2t-s\sep t+1\dots 2t].
$$
We then apply relation in Lemma \ref{genplu2}. An easy check shows
that we get an expression of $\delta_tM$ as a sum of terms $\pm
N_1N_2$ such that $N_1$ and $N_2$ are $t$-minors with $e(N_i)>e(M)$
(coordinate wise). By induction, we conclude that there exists an
exponent $h$ such that $f^hN_1N_2\in K[\Phi_2]$. So
$\delta_t^{h+1}M\in K[\Phi_2]$.

The cases (a) and (b) are illustrated by the relations
\eqref{12terms3min} and \eqref{12terms2min} where the underlined
terms correspond to $\delta_tM$.\medskip

(2)  Let $M$ be a $t$-minor in $\Phi_2\setminus \Phi_1$. Then $M$
contains at least $t-1$ row indices $\leq t$ and hence at least
$t-2$ row indices $\leq t-1$, and the same holds for columns. The
statement we have to prove is completely symmetric in the first
$t-1$ row and column indices. So we may assume that the row indices
of $M$ are $1,2,\dots, t-2,r_1,r_2$ and the column indices are
$1,2,\dots, t-2,c_1,c_2$. By the retraction principle, we may assume
that $t=2$. Then, up to transposition and renaming of the indices
larger than $3$, $M$ is one of the following:
$$
[13|24], \quad [14|23], \quad [14|24], \quad [23|24], \quad [24|24].
$$
The relations \eqref{6a}--\eqref{6e} indeed imply that
$\delta_1\delta_3$ times each element of the list above belongs to
the $K$-algebra generated by $\Phi_1$ and elements which are earlier
in the list (up to transposition).

We verify this in detail: relation \eqref{6e} yields that
$\delta_1\delta_3[13|24]$ belongs to $K[\Phi_1]$. Relation
\eqref{6d} says that $\delta_1\delta_3[14|23]$ belongs to
$K[\Phi_1]$. Relation \eqref{6c} asserts that
$\delta_1\delta_3[14|24]$ is in $K\bigl[\Phi_1, [13|24], [14|23]
\bigr]$. Relation \eqref{6b} means that $\delta_1\delta_3[23|24]$
belongs to $K\bigl[\Phi_1,\allowbreak [13|24],\allowbreak
[23|14]\bigr]$. Relation \eqref{6a} guarantees that
$\delta_1\delta_3[24|24]$ lies in $K\bigl[\Phi_1, [13|24],
[24|13],\allowbreak [14|23],\allowbreak [23|14]\bigr]$. This
concludes the proof of Claim (2).\medskip

(3) Let $M\in \Phi_1\setminus \Phi_0$. Up to transposition,
$M=[1\dots t\sep b_1b_2\dots b_t]$ with $b_{t-1}>t$. Set $e(M)=\{ i
: b_i>t\}$. Then $e(M)>1$. Using  the Pl\"ucker relations one can
write $\delta_tM$ as a linear combination of products $N_1N_2$ with
$N_i=[1\dots t\sep c_1c_2\dots c_t]$ and $N_2=[1\dots t\sep
c_1c_2\dots c_t]$ with $e(N_1)<e(M)$ and $e(N_2)<e(M)$. Iterating
the arguments, one concludes that $\delta_t^kM\in K[\Phi_0]$ for
some $k$. The details are given in \cite[6.1]{BV}.
\end{proof}

The proof of the Theorem \ref{localize} is now complete. Our next
goal is to identify a singular points of the variety $X_t$. We will
use the following:

\begin{lemma}\label{nonsingular}
Let $R$ be a regular ring and $I,\pp$ ideals of $R$ such that
$I\subset \pp$ and $\pp$ is prime. Set $N=(I+\pp^2)/\pp^2$. If
$R_\pp/IR_\pp$ is regular, then the $R$-module $N$ is minimally
generated by at least $\height I$ elements.
\end{lemma}

\begin{proof} Since $R_\pp/IR_\pp$ is regular, the ideal $IR_\pp$ is
minimally generated by $\height I$ elements belonging to a regular
system of parameters of $R_\pp$ (for example, see \cite[2.2.4]{BH}).
In other words, $N_\pp=(I+\pp^2)R_\pp/\pp^2R_\pp$ is generated, as
an $R_\pp$-module, by $\height I$ elements. Since the number of
generators can only decrease under localization, $N$ is minimally
generated by at least $\height I$ elements.
\end{proof}

Now we can prove that $X_t$ is singular in certain points.

\begin{proposition}\label{non-iso-sing}
Suppose that $1<t<min(m,n)$, and that $t\neq m-1$ if $m=n$, and let
$x\in X_t(m,n)$ be a point with $\rank x=1$. Then $X_t(m,n)$ is
singular at $x$.
\end{proposition}

\begin{proof}
Set $A=A_t$. We may assume $x=\phi(\alpha)$ where $\alpha$ is a
$m\times n$ matrix whose only non-zero $t$-minor is
$\delta=[1\dots,t\sep 1\dots t]$. Then the maximal ideal
corresponding to $x$ contains the ideal $\pp$ of $A$ generated by
all the $t$-minors $\gamma\neq \delta$. It is enough to prove that
$A_\pp$ is not regular. For each $t$-minor $\gamma$ of $X$ we pick a
new variable $Y_\gamma$ and present $A$ in the form $R/I$ with
$R=K[Y_\gamma: \gamma\text{ $t$-minor of $X$}]$. We give $R$ the
standard $\ZZ^m\times \ZZ^n$-grading. Let $\pp$ denote the ideal of
$R$ generated by all $Y_\gamma$ with $\gamma\neq \delta$. According
to Lemma \ref{nonsingular} it is enough to show that
$N=(I+\pp^2)/\pp^2$ is generated by fewer than $\height I$ elements
as an $R$-module. Note that
$$
\height I=\binom{m}{t}\binom{n}{t}-mn.
$$
Note also that each non-zero $\ZZ^m\times \ZZ^n$-homogeneous element
in $R/\pp^2$ is, up to irrelevant scalars, of the form
$Y_\delta^kY_\gamma \mod \pp^2$. Since $I$ is $\ZZ^m\times \ZZ^n$
homogeneous, it follows that a $K$-basis of $N$ is given by the
elements $Y_\delta^kY_\gamma \mod \pp^2$ such that there exist $g\in
I$ of the form $g= Y_\delta^kY_\gamma+f$ with $f\in \pp^2$. Then the
generators of $N$ as an $R$-module are the elements
$Y_\delta^kY_\gamma \mod \pp^2$ such that $g$ is as above and $k$ is
taken minimal for the given $\gamma$. Finally note that if $\gamma$
and $\delta$ are contained in a $(t+1)\times (t+1)$ submatrix of $X$
the element $Y_\delta^kY_\gamma \mod \pp^2$ cannot be in $N$ because
the $t$-minors of a $(t+1)\times (t+1)$ matrix of indeterminates are
algebraically independent. The number of the $t$-minors $\gamma$
such that $\gamma$ and $\delta$ are contained in a $(t+1)\times
(t+1)$ submatrix is $(t(m-t)+1)(t(n-t)+1)$. So the number of
generators of $N$ is at most
$$
\binom{m}{t}\binom{n}{t}-(t(m-t)+1)(t(n-t)+1).
$$
Therefore we have to show that
$$
(t(m-t)+1)(t(n-t)+1)>mn.
$$
To this end, it suffices that
$$
t(m-t)+1\geq m \quad\text{and}\quad t(n-t)+1>n
$$
which is equivalent to $m\geq t+1$ and $n>t+1$, and so we are done.
\end{proof}

Before we summarize the results of this section, let us discuss
those cases that have been excluded in Proposition
\ref{non-iso-sing}. First, if $t=1$ or $t=m-1=n-1$, then $X_t(m,n)$
is the full affine space. This is trivial for $t=1$, and for
$t=m-1=n-1$ we have already seen this in the introduction. If $t=m$,
then $X_t(m,n)$ is the cone over the Grassmannian, and so $0$ is the
only singular point of $X_t(m,n)$. In the remaining cases the
singular locus is given by the next theorem.

\begin{theorem}\label{SingLoc}
Suppose that $1<t<min(m,n)$, and that $t\neq m-1$ if $m=n$. Then the
singular locus of $X_t(m,n)$ consists of the points $x$ such that
$\rank x\le 1$.
\end{theorem}

\begin{proof}
As we have seen in Proposition \ref{non-iso-sing}, $X=X_t(m,n)$ is
indeed singular at all points of rank $1$ and (therefore) also at
$x=0$.

The singular locus of $X$ is $G$-stable. Therefore its components
are defined by $G$-stable prime ideals $\pp$ for which $(A_t)_\pp$
is non-regular. The prime ideal $\qq=\qq_{t+1}$ defines the locus of
points of rank $\le1$ (Corollary \ref{all-pairs}). Therefore it is
enough to show that $(A_t)_\pp$ is regular for all $G$-stable prime
ideals $\pp\neq \qq_t,\qq_{t+1}$.  As a consequence of Lemma
\ref{primes3} such $\pp$ does not contain the denominator $F$ of
Theorem \ref{localize}, and it follows immediately that $(A_t)_\pp$
is regular.
\end{proof}

\begin{remark}
As we have seen, up to the singular locus, $X_t$ is defined by
equations of degree $\le 3$. It seems that $X_t$ itself is defined
by such equations. At least we do not know of a counterexample.
However, equations of degree $2$ are not sufficient if $1<t<m$ and
$t\neq m-1=n-1$, as is demonstrated by the case $m=3$, $n=4$, $t=2$
discussed in \cite{B}.
\end{remark}

\begin{remark}
Theorem \ref{localize} and Proposition \ref{non-iso-sing} hold in
arbitrary characteristic.
\end{remark}


\begin{thebibliography}{99}

\bibitem{Bou:alg} N. Bourbaki, {\em Alg\`ebre, Ch.\ 1 \`a 3}.
Hermann 1970.

\bibitem{B} W. Bruns, {\em Algebras defined by powers of
determinantal ideals}, J. of Algebra 142 (1991), 150--163.

\bibitem{BC2} W. Bruns and A. Conca, {\em Algebras of minors}.
 J. Algebra 246 (2001), 311--330.

 \bibitem{BC3} W. Bruns and A. Conca. {\em
The F-rationality of determinantal rings and their Rees rings.}
Mich. Math. J. {\bf 45} (1998), 291--299.

\bibitem{BH}  W. Bruns and  J. Herzog, {\em Cohen-Macaulay rings},
Cambridge Studies in Advanced Mathematics, 39,  Cambridge University
Press 1993.

\bibitem{BV} W. Bruns and U. Vetter, {\em Determinantal rings},
Lect. Notes Math. 1327, Springer 1988.

\bibitem{DEP} C. De Concini, D. Eisenbud, and C. Procesi, {\em Young
diagrams and determinantal varieties}, Invent. math. 56 (1980),
129--165.

\bibitem{DRS} P. Doubilet, G.C. Rota, and J. Stein, {\em On the
foundations of combinatorial theory: IX, Combinatorials methods in
invariants theory}, Studies in Applied Mathematics LIII (1974),
185--216.

\bibitem{Gr} F. D. Grosshans, {\em Contractions of the actions of reductive
algebraic groups in arbitrary charactersitic}, Invent. math. {\bf
107} (1992), 127--133.

\bibitem{St} R. Steinberg, {\em Conjugacy classes in algebraic groups},
Lecture Notes in Mathematics 366. Springer 1974)



\end{thebibliography}
\end{document}